\documentclass[11pt]{article}

\usepackage[margin=1in]{geometry}
\usepackage{amsmath,amssymb,amsthm,mathtools,mathrsfs}
\usepackage{booktabs}
\usepackage{graphicx}
\usepackage{microtype}
\usepackage[hidelinks]{hyperref}

\newtheorem{theorem}{Theorem}[section]
\newtheorem{proposition}[theorem]{Proposition}
\newtheorem{lemma}[theorem]{Lemma}
\newtheorem{corollary}[theorem]{Corollary}
\newtheorem{definition}[theorem]{Definition}
\newtheorem{remark}[theorem]{Remark}

\newcommand{\R}{\mathbb R}
\newcommand{\N}{\mathbb N}

\newcommand{\T}{\mathcal T}
\newcommand{\A}{\mathcal A}
\newcommand{\Bf}{\mathcal B}
\newcommand{\Lap}{\mathscr L}
\newcommand{\Dt}{\mathscr D}
\newcommand{\G}{G}
\newcommand{\bc}{\mathbf c}
\newcommand{\bn}{\boldsymbol\nu}
\newcommand{\bm}{\boldsymbol\mu}

\title{Generalized Images and Multi-Index Transmutations\\
for Brownian First Passage to Polynomial Boundaries}

\author{Gerardo Hern\'andez-del-Valle\thanks{Centro de Estudios Monetarios
Latinoamericanos (CEMLA), Mexico City, Mexico. Corresponding author.
E-mail: \texttt{ghernandez@cemla.org}}}
\date{}

\begin{document}
\maketitle

\begin{abstract}
We study Brownian first passage to a finite polynomial boundary
\[
 F(t)=d+\sum_{j=2}^p c_jt^j.
\]
For every such boundary we construct a time-independent signed exterior
image measure.  After Brownian scaling its density is determined by a
second-kind Volterra equation whose kernel is the multivariate Wright series
\[
 \mathcal K_{\boldsymbol\beta}(u,v)
 =\sum_{\bn\ne\mathbf0}
 \frac{\boldsymbol\beta^{\bn}v^{|\bn|/2}
 (u-v)^{\kappa(\bn)-1}}
 {2^{\kappa(\bn)}\bn!\Gamma(\kappa(\bn))},
 \qquad
 \kappa(\bn)=\sum_{j=2}^p(j-1)\nu_j.
\]
An exact multivariate finite-difference cancellation yields exponential tail
control of the image density.  Consequently, on a nonzero interval of the
natural Brownian scale, its heat potential is the killed transition kernel
and its boundary flux is the first-passage density.

The image kernel is connected exactly with the differential transported-zero
theory.  A two-variable beta-moment transform carries its multivariate Wright
seed to the entire constant-flux function associated with the polynomial
boundary.  Writing $K(\bn)=\sum j\nu_j$ and
$W(\bn)=\sum(2j-1)\nu_j$, we construct the graded spatial transmutation and
prove the Laplace conjugacy
\[
 \Lap\A_{\bn}\Lap^{-1}
 =\frac{(-1)^{|\bn|}}{\bn!}k^{|\bn|}
 \left(-\frac1k\frac{d}{dk}\right)^{K(\bn)}.
\]
Thus the inverse transmutation generates the formal parameter expansion of
the same first-passage problem.  For $F(t)=At^2+Bt^3$ we also derive
fifth- and eighth-order scalar annihilators and prove their exact rational
factorization.  Numerical experiments for two- and three-term polynomial
boundaries independently agree with Durbin's equation to discretization
accuracy.
\end{abstract}



\section{Introduction}

Let $(B_t)_{t\ge0}$ be standard Brownian motion started at $B_0=-a$, where
$a>0$, and consider the finite polynomial boundary
\begin{equation}\label{eq:F}
 F_{\bc}(t)=\sum_{j=2}^p c_jt^j,
 \qquad \bc=(c_2,\ldots,c_p)\in\R^{p-1}.
\end{equation}
The corresponding first-passage time is
\begin{equation}\label{eq:tau}
 \tau_{\bc,a}=\inf\{t>0:B_t\ge F_{\bc}(t)\}.
\end{equation}
For a general curved boundary, the law of \eqref{eq:tau} is described by
integral equations, moving-domain heat problems, changes of measure and
approximation methods; see, among others, \cite{Lerche,Durbin,Peskir,JKV,
HerrmannTanre}.  Explicit special-function descriptions occur only for
restricted families.  The quadratic case, for example, is closely connected
with Airy functions \cite{Salminen,Groeneboom}.  Boundary transformations
and nonlinear ordinary differential equations provide another important
perspective \cite{AliliPatie}.

The inverse method of images seeks a time-independent measure whose heat
potential enforces the prescribed moving boundary.  Analytic examples and a
recent general duality framework are given in \cite{Kahale,Christensen}; the
existence of first-passage densities for sufficiently regular curves is
treated, for example, in \cite{Lee}.  Positivity of the image measure is a
substantive restriction.  The polynomial construction below instead uses a
signed exterior measure.

The present paper approaches \eqref{eq:tau} through the transported-zero
problem for the heat semigroup
\[
 P_t=\exp\!\left(\frac t2\partial_x^2\right).
\]
Given a curve $F$, one seeks an initial function $\Gamma$ such that
\begin{equation}\label{eq:moving-zero}
 (P_t\Gamma)(F(t))=0.
\end{equation}
This problem arose from attempts to construct heat solutions adapted to
Brownian boundary crossing \cite{HDVGuerra}.  Transported differential
operators, compatibility matrices and auxiliary Heat--Airy flows organize
the direct problem of determining $F$ from a stationary equation
\cite{HDVMatrices,HDVHeatAiry}.  Stationary grading explains why a naive
order-$p$ stationary ansatz produces exact quadratic and cubic zeros but
fails at the quartic level \cite{HDVGrading}.  The inverse problem instead
prescribes $F$ and reconstructs admissible initial functions
\cite{HDVInverse}.

For the monomial boundary $F(t)=\rho t^p$, the constant-flux normalization
\begin{equation}\label{eq:flux-normalization-intro}
 (P_t\Gamma)(\rho t^p)=0,
 \qquad
 \partial_x(P_t\Gamma)(\rho t^p)=1
\end{equation}
selects a hypergeometric function annihilated by a polynomial-coefficient
differential operator of order $2p-1$ \cite{HDVMonomial}.  The associated
graded spatial transmutation is Laplace-conjugate to the temporal moment
operator in the Brownian exponential-martingale identity.  Its inverse
therefore generates the complete formal first-passage expansion for the
monomial boundary \cite{HDVTransmutation}.

The probabilistic normalization requires something different.  A killed
Brownian kernel must begin from an interior Dirac mass, whereas the
constant-flux selection is function-valued.  For monomial boundaries this
gap was closed by constructing a time-independent signed exterior image
measure through a second-kind Volterra equation in the image coordinate
\cite{HDVImages}.  Its Wright-function kernel is linked by beta moments to
the monomial constant-flux function.  The present paper shows that both the
generalized-image construction and its differential structure extend from a
single monomial to every finite polynomial.

The first main result is probabilistic.  For
\[
 F(t)=d+\sum_{j=2}^pc_jt^j,
 \qquad B_0=d-a,
\]
we construct a signed measure supported in $[d+a,\infty)$ whose heat potential
vanishes on $F$.  After scaling, its density solves an explicit multivariate
Volterra equation.  A mixed finite-difference cancellation extends the
factorization method of the monomial case and gives exponential control at
the rate
\[
 \max_{2\le j\le p}|c_ja^{2j-1}|^{2/(2j-1)}.
\]
The resulting heat potential is therefore the killed Brownian transition
kernel on a nonzero interval, and its boundary flux is the actual
first-passage density.

The second main result identifies the differential structure hidden in this
image construction.  Every finite polynomial \eqref{eq:F} determines an
explicit multi-index transmutation.  Its grades are governed by
\[
 K(\bn)=\sum_{j=2}^p j\nu_j,
 \qquad
 W(\bn)=\sum_{j=2}^p(2j-1)\nu_j.
\]
Here $K(\bn)$ is the temporal moment order, whereas $W(\bn)$ is the spatial
degree shift.  Their relation
\[
 W(\bn)=2K(\bn)-|\bn|
\]
is precisely the Brownian scaling balance required by Laplace conjugacy.

At the kernel level, introduce
\[
 \kappa(\bn)=\sum_{j=2}^p(j-1)\nu_j.
\]
Then $W(\bn)=K(\bn)+\kappa(\bn)$.  The image-coordinate Wright seed has
coefficients $1/[\bn!\Gamma(\kappa(\bn))]$.  A multivariate beta-moment
operator multiplies these coefficients by
\[
 \frac{\Gamma(K(\bn)+1)\Gamma(\kappa(\bn))}
 {\Gamma(W(\bn)+2)},
\]
and hence produces exactly the constant-flux function.  The generalized
image theory and the transported-zero theory therefore share an exact
transform, not only a common scaling exponent.

Three features distinguish the polynomial theory from a formal superposition
of the monomial cases.  First, the transmutation is indexed by the full
multi-index $\bn$.  Second, the mixed operators are not compositions of the
one-variable operators.  Third, the inverse recursion contains genuine
interactions between different boundary powers.  Already for
$F(t)=At^2+Bt^3$, the coefficient of $AB$ materially improves the numerical
first-passage approximation.

The remaining results organize this connection.  The transmutation transports
every odd flat-boundary heat mode to $F_{\bc}$ while preserving its boundary
flux.  Its image of $x$ is the entire constant-flux function obtained from
the Wright seed.  Exact Laplace conjugacy identifies its inverse spatial
recursion with the formal recursion from the exponential martingale.  In the
quadratic--cubic case we derive scalar ODEs of orders five and eight.  Their
factorization shows why differential order and perturbative usefulness are
distinct.  Finally, numerical image calculations for quadratic--cubic
and quadratic--cubic--quartic boundaries agree independently with Durbin's
equation, while the transmutation expansion resolves the contribution of the
mixed coefficients.

\section{Polynomial exterior images}
\label{sec:images}

Write
\[
 f_{\bc}(t)=\sum_{j=2}^pc_jt^j,
 \qquad
 F(t)=d+f_{\bc}(t),
 \qquad
 x_0=d-a,
 \qquad a>0,
\]
and let
\[
 \G(t,x)=\frac1{\sqrt{2\pi t}}\exp\!\left(-\frac{x^2}{2t}\right).
\]
We seek a time-independent signed measure $\mu_{\bc,a}$ supported to the
right of the initial domain such that
\begin{equation}\label{eq:image-potential}
 P_{\bc,a}(t,x)
 =\G(t,x-x_0)-\int_d^\infty\G(t,x-y)\,\mu_{\bc,a}(dy)
\end{equation}
satisfies $P_{\bc,a}(t,F(t))=0$.

Put $y=d+z$.  Cancelling the common factor generated by
$f_{\bc}(t)^2$ in the boundary identity gives
\begin{equation}\label{eq:polynomial-image-unscaled}
 \exp\!\left(-\frac{a^2}{2t}-\frac{af_{\bc}(t)}t\right)
 =\int_0^\infty
 \exp\!\left(-\frac{z^2}{2t}+\frac{zf_{\bc}(t)}t\right)
 \,\widetilde\mu_{\bc,a}(dz).
\end{equation}
Scale
\begin{equation}\label{eq:polynomial-image-scaling}
 z=ar,
 \qquad t=a^2s,
 \qquad \lambda=\frac1{2s},
 \qquad \beta_j=c_ja^{2j-1},
\end{equation}
and push the measure forward under $r\mapsto u=r^2$.  The causal Laplace
inversion below shows that the resulting measure is supported on
$[1,\infty)$.  Equation~\eqref{eq:polynomial-image-unscaled} becomes
\begin{equation}\label{eq:polynomial-image-transform}
 \exp\!\left(
 -\lambda-
 \sum_{j=2}^p\frac{\beta_j}{2^{j-1}\lambda^{j-1}}
 \right)
 =\int_1^\infty
 \exp\!\left(
 -\lambda u+\sqrt u
 \sum_{j=2}^p\frac{\beta_j}{2^{j-1}\lambda^{j-1}}
 \right)\nu_{\boldsymbol\beta}(du).
\end{equation}

For $\bn=(\nu_2,\ldots,\nu_p)$ define
\begin{equation}\label{eq:kappa-image}
 |\bn|=\sum_{j=2}^p\nu_j,
 \qquad
 \kappa(\bn)=\sum_{j=2}^p(j-1)\nu_j.
\end{equation}
For $u>v\ge1$, put
\begin{equation}\label{eq:polynomial-image-kernel}
 \boxed{
 \mathcal K_{\boldsymbol\beta}(u,v)
 =\sum_{\bn\ne\mathbf0}
 \frac{\boldsymbol\beta^{\bn}v^{|\bn|/2}
 (u-v)^{\kappa(\bn)-1}}
 {2^{\kappa(\bn)}\bn!\Gamma(\kappa(\bn))}.}
\end{equation}
The series is locally uniformly convergent on the Volterra triangle.  Let
\[
 (\mathcal K_{\boldsymbol\beta}h)(u)
 =\int_1^u\mathcal K_{\boldsymbol\beta}(u,v)h(v)\,dv,
 \qquad
 \mathcal T^{\rm im}_{\boldsymbol\beta}
 =I+\mathcal K_{\boldsymbol\beta}.
\]

Indeed, on a compact triangle $1\le v\le u\le U$, the term indexed by
$\bn$ is bounded by
\[
 \frac1{U\bn!\Gamma(\kappa(\bn))}
 \prod_{j=2}^p
 \left(\frac{|\beta_j|\sqrt U\,U^{j-1}}{2^{j-1}}\right)^{\nu_j}.
\]
Since $\kappa(\bn)\ge|\bn|$, this has a summable majorant.  In particular,
the apparent factor $(u-v)^{-1}$ in \eqref{eq:kernel-Wright} is removable;
the original series \eqref{eq:polynomial-image-kernel} is the definition on
the diagonal.

\begin{proposition}[Polynomial image-coordinate equation]
\label{prop:polynomial-image-equation}
Equation~\eqref{eq:polynomial-image-transform} is equivalent, in the class
of locally finite measures of suitable exponential type, to
\begin{equation}\label{eq:polynomial-volterra-measure}
 \boxed{
 \mathcal T^{\rm im}_{\boldsymbol\beta}\nu_{\boldsymbol\beta}
 =\mathcal T^{\rm im}_{-\boldsymbol\beta}\delta_1.}
\end{equation}
Writing
\[
 \nu_{\boldsymbol\beta}
 =\delta_1+n_{\boldsymbol\beta}(u)\,du,
\]
the continuous image density solves
\begin{equation}\label{eq:polynomial-volterra-density}
 \boxed{
 n_{\boldsymbol\beta}(u)
 +\int_1^u\mathcal K_{\boldsymbol\beta}(u,v)
 n_{\boldsymbol\beta}(v)\,dv
 =\mathcal K_{-\boldsymbol\beta}(u,1)
 -\mathcal K_{\boldsymbol\beta}(u,1).}
\end{equation}
\end{proposition}

\begin{proof}
Expand the second exponential in
\eqref{eq:polynomial-image-transform}.  For every nonzero multi-index $\bn$,
\[
 \mathcal L^{-1}_{\lambda\to u}
 \left(e^{-\lambda v}\lambda^{-\kappa(\bn)}\right)
 =\frac{(u-v)_+^{\kappa(\bn)-1}}
 {\Gamma(\kappa(\bn))}.
\]
Termwise inversion gives \eqref{eq:polynomial-volterra-measure}; separating
the atom at $u=1$ gives \eqref{eq:polynomial-volterra-density}.  The local
uniform kernel bound above justifies the inversion first on compact
Volterra triangles.  For measures of exponential type, both Laplace
transforms converge for all sufficiently large $\lambda$, where uniqueness
of the Laplace transform proves the converse.
\end{proof}

\section{The multivariate Wright seed and beta moments}
\label{sec:beta-moments}

The kernel in \eqref{eq:polynomial-image-kernel} is generated by the
following multivariate extension of Wright's generalized Bessel series
\cite{Wright}:
\begin{equation}\label{eq:multi-Wright}
 \mathcal W_p(\mathbf z)
 =\sum_{\bn\ne\mathbf0}
 \frac{\mathbf z^{\bn}}
 {\bn!\Gamma(\kappa(\bn))}.
\end{equation}
Indeed,
\begin{equation}\label{eq:kernel-Wright}
 \mathcal K_{\boldsymbol\beta}(u,v)
 =\frac1{u-v}\mathcal W_p\!\left(
 \left(
 \frac{\beta_j\sqrt v\,(u-v)^{j-1}}{2^{j-1}}
 \right)_{j=2}^p
 \right).
\end{equation}

For a power series $\Phi$ with $\Phi(\mathbf0)=0$, define
\begin{equation}\label{eq:multivariate-beta-operator}
 (\mathfrak B_p\Phi)(\mathbf z)
 =\int_0^1\!\int_0^1
 \frac{
 \Phi((z_jt^j(1-t)^{j-1}s^{2j-1})_{j=2}^p)}{1-t}
 \,ds\,dt.
\end{equation}

\begin{proposition}[Multivariate Wright-to-constant-flux transform]
\label{prop:multi-beta-transform}
Let $K(\bn)$ and $W(\bn)$ be as in \eqref{eq:KW}.  Then
\begin{equation}\label{eq:multi-beta-series}
 1+\mathfrak B_p\mathcal W_p(\mathbf z)
 =\sum_{\bn\ge\mathbf0}
 \frac{\Gamma(K(\bn)+1)}
 {\bn!\Gamma(W(\bn)+2)}\mathbf z^{\bn}.
\end{equation}
Consequently,
\begin{equation}\label{eq:Gamma-beta-transform}
 \boxed{
 \Gamma_{\bc}(x)
 =x\left[1+\mathfrak B_p\mathcal W_p
 \bigl((-2^jc_jx^{2j-1})_{j=2}^p\bigr)\right]}
\end{equation}
is exactly the constant-flux function constructed in
Corollary~\ref{cor:Gamma} below.
\end{proposition}

\begin{proof}
The defining series of $\mathcal W_p$ is normally convergent on compact
subsets because $\kappa(\bn)\ge|\bn|$.  Since it has no constant term, the
factor $(1-t)^{-1}$ in \eqref{eq:multivariate-beta-operator} is cancelled by
at least one power of $(1-t)$, so termwise integration is justified on
compact $\mathbf z$-sets.
On the monomial $\mathbf z^{\bn}$, the two integrations in
\eqref{eq:multivariate-beta-operator} give
\begin{align*}
 &\int_0^1t^{K(\bn)}(1-t)^{\kappa(\bn)-1}\,dt
 \int_0^1s^{W(\bn)}\,ds\\
 &\hspace{20mm}=
 \frac{\Gamma(K(\bn)+1)\Gamma(\kappa(\bn))}
 {\Gamma(W(\bn)+2)},
\end{align*}
because $W(\bn)=K(\bn)+\kappa(\bn)$.  Multiplication by the coefficient in
\eqref{eq:multi-Wright} proves \eqref{eq:multi-beta-series}.  Substituting
$z_j=-2^jc_jx^{2j-1}$ yields \eqref{eq:Gamma-beta-transform}.
\end{proof}

\section{Factorization and multivariate cancellation}
\label{sec:image-tail}

Apply $\mathcal T^{\rm im}_{-\boldsymbol\beta}$ to
\eqref{eq:polynomial-volterra-measure} and write
\begin{equation}\label{eq:image-factorization}
 \mathcal T^{\rm im}_{-\boldsymbol\beta}
 \mathcal T^{\rm im}_{\boldsymbol\beta}
 =I+\mathcal E_{\boldsymbol\beta}.
\end{equation}
The factorized equation is
\begin{equation}\label{eq:image-factorized-equation}
 (I+\mathcal E_{\boldsymbol\beta})\nu_{\boldsymbol\beta}
 =(\mathcal T^{\rm im}_{-\boldsymbol\beta})^2\delta_1.
\end{equation}

\begin{proposition}[Exact multivariate defect]
\label{prop:multi-defect}
The defect begins at total parameter degree two.  For
$|\boldsymbol\alpha|\ge2$, the coefficient of
$\boldsymbol\beta^{\boldsymbol\alpha}$ is
\begin{align}
 E_{\boldsymbol\alpha}(u,v)
 ={}&\frac1{2^{\kappa(\boldsymbol\alpha)}}
 \sum_{\mathbf0<\boldsymbol\eta<\boldsymbol\alpha}
 \frac{(-1)^{|\boldsymbol\eta|}
 v^{|\boldsymbol\alpha-\boldsymbol\eta|/2}}
 {\boldsymbol\eta!(\boldsymbol\alpha-\boldsymbol\eta)!
 \Gamma(\kappa(\boldsymbol\eta))
 \Gamma(\kappa(\boldsymbol\alpha-\boldsymbol\eta))}
 \notag\\
 &\quad\times\int_v^u
 \left(w^{|\boldsymbol\eta|/2}
 -v^{|\boldsymbol\eta|/2}\right)
 (u-w)^{\kappa(\boldsymbol\eta)-1}
 (w-v)^{\kappa(\boldsymbol\alpha-\boldsymbol\eta)-1}
 \,dw.
 \label{eq:multi-defect-coefficient}
\end{align}
\end{proposition}

\begin{proof}
Expand
$\mathcal K_{-\boldsymbol\beta}
+\mathcal K_{\boldsymbol\beta}
+\mathcal K_{-\boldsymbol\beta}
 \mathcal K_{\boldsymbol\beta}$ by multi-index.  In the composition term,
replace $w^{|\boldsymbol\eta|/2}$ temporarily by
$v^{|\boldsymbol\eta|/2}$.  The beta integral then produces the same kernel
monomial as the two single-kernel terms.  Its coefficient cancels because
\[
 \sum_{\boldsymbol\eta\le\boldsymbol\alpha}
 (-1)^{|\boldsymbol\eta|}
 \binom{\boldsymbol\alpha}{\boldsymbol\eta}=0.
\]
Restoring the difference gives \eqref{eq:multi-defect-coefficient}.
\end{proof}

\begin{lemma}[Mixed finite-difference cancellation]
\label{lem:mixed-cancellation}
For $m\ge1$, define
\begin{equation}\label{eq:mixed-cancellation-sum}
 S_{\boldsymbol\alpha,m}
 =\sum_{\mathbf0<\boldsymbol\eta<\boldsymbol\alpha}
 \frac{(-1)^{|\boldsymbol\eta|}}
 {\boldsymbol\eta!(\boldsymbol\alpha-\boldsymbol\eta)!}
 \binom{|\boldsymbol\eta|/2}{m}
 (\kappa(\boldsymbol\alpha-\boldsymbol\eta))_m,
\end{equation}
where $(z)_m=z(z+1)\cdots(z+m-1)$.  Then
\begin{equation}\label{eq:mixed-cancellation-zero}
 \boxed{S_{\boldsymbol\alpha,m}=0
 \qquad\text{whenever}\qquad 2m<|\boldsymbol\alpha|.}
\end{equation}
\end{lemma}

\begin{proof}
After multiplication by $\boldsymbol\alpha!$, the summand is the restriction
to the integer box $\mathbf0\le\boldsymbol\eta\le\boldsymbol\alpha$ of a
polynomial in $\boldsymbol\eta$ of total degree at most $2m$.  The two
endpoint terms vanish and may be included.  The resulting sum is the mixed
finite difference of that polynomial of total order
$|\boldsymbol\alpha|$.  It is zero when $2m<|\boldsymbol\alpha|$.
\end{proof}

Put
\[
 m_n=\left\lceil\frac n2\right\rceil,
 \qquad
 b_{n,p}=\frac{e^{m_n}(pn)^{m_n}}{n!}.
\]
Stirling's formula gives
\begin{equation}\label{eq:bn-superexponential}
 \limsup_{n\to\infty}b_{n,p}^{1/n}=0.
\end{equation}

\begin{proposition}[Global defect estimate]
\label{prop:global-defect}
For $n=|\boldsymbol\alpha|\ge2$,
\begin{align}
 |E_{\boldsymbol\alpha}(u,v)|
 \le{}&
 \frac{e^{m_n}(pn)^{m_n}
 2^{n-\kappa(\boldsymbol\alpha)}}
 {\boldsymbol\alpha!}
 v^{n/2-m_n}
 \notag\\
 &\times
 \frac{(u-v)^{\kappa(\boldsymbol\alpha)+m_n-1}}
 {\Gamma(\kappa(\boldsymbol\alpha)+m_n)}.
 \label{eq:global-defect-estimate}
\end{align}
\end{proposition}

\begin{proof}
Put $r=|\boldsymbol\eta|/2$.  Taylor's formula gives
\[
 w^r=\sum_{m=0}^{m_n-1}
 \binom rm v^{r-m}(w-v)^m+R_{\boldsymbol\eta}(w,v).
\]
For every $m<m_n$ one has $2m<n$, so after the beta integration the
coefficient of the $m$th Taylor term is precisely
$S_{\boldsymbol\alpha,m}$ and vanishes by
Lemma~\ref{lem:mixed-cancellation}; the case $m=0$ is the binomial
cancellation already used in Proposition~\ref{prop:multi-defect}.  Since
$0\le r\le n/2\le m_n$ and $v\ge1$, the integral form of the remainder
implies the uniform bound
\[
 |R_{\boldsymbol\eta}(w,v)|
 \le e^{m_n}v^{r-m_n}(w-v)^{m_n}.
\]
Inserting this estimate into \eqref{eq:multi-defect-coefficient}, the beta
integral contributes
$(\kappa(\boldsymbol\alpha-\boldsymbol\eta))_{m_n}$ and the common factor
\[
 \frac{(u-v)^{\kappa(\boldsymbol\alpha)+m_n-1}}
 {\Gamma(\kappa(\boldsymbol\alpha)+m_n)}.
\]
Moreover,
$\kappa(\boldsymbol\alpha-\boldsymbol\eta)+m_n-1\le pn$, so the rising
factorial is at most $(pn)^{m_n}$.  Finally,
\[
 \sum_{\boldsymbol\eta\le\boldsymbol\alpha}
 \frac1{\boldsymbol\eta!(\boldsymbol\alpha-\boldsymbol\eta)!}
 =\frac{2^n}{\boldsymbol\alpha!}.
\]
This proves \eqref{eq:global-defect-estimate}.
\end{proof}

For $\sigma>0$ and $V\ge1$, let
\[
 X_{\sigma,V}
 =\left\{h:
 \|h\|_{\sigma,V}
 =\sup_{u\ge V}e^{-\sigma(u-V)}|h(u)|<\infty\right\}.
\]
Define
\[
 z_j=\frac{|\beta_j|}{\sigma^{j-1/2}},
 \qquad Z=\sum_{j=2}^pz_j,
\]
and the entire functions
\[
 \Phi_p(Z)=\sum_{\substack{n\ge2\\n\ \mathrm{even}}}b_{n,p}Z^n,
 \qquad
 \Psi_p(Z)=\sum_{\substack{n\ge3\\n\ \mathrm{odd}}}b_{n,p}Z^n.
\]
Here we used the multinomial identity
\[
 \sum_{|\boldsymbol\alpha|=n}
 \frac{\mathbf z^{\boldsymbol\alpha}}{\boldsymbol\alpha!}
 =\frac{Z^n}{n!}.
\]
To obtain the operator norm, multiply the kernel estimate by
$e^{-\sigma(u-v)}$ and use
\[
 \int_0^\infty e^{-\sigma s}
 \frac{s^{q-1}}{\Gamma(q)}\,ds=\sigma^{-q}.
\]
Since
\[
 \kappa(\boldsymbol\alpha)+\frac n2
 =\sum_{j=2}^p\left(j-\frac12\right)\alpha_j,
\]
the even levels give $\mathbf z^{\boldsymbol\alpha}$, while the odd levels
give the additional factor $(\sigma v)^{-1/2}\le(\sigma V)^{-1/2}$.
The parity in \eqref{eq:global-defect-estimate} therefore gives
\begin{equation}\label{eq:defect-tail-norm}
 \|\mathcal E_{\boldsymbol\beta}^{(V)}\|_{\sigma,V}
 \le\Phi_p(Z)+(\sigma V)^{-1/2}\Psi_p(Z).
\end{equation}

\begin{theorem}[Polynomial image density]
\label{thm:polynomial-image-density}
For every $\boldsymbol\beta\in\R^{p-1}$,
\eqref{eq:polynomial-volterra-measure} has a unique locally finite solution
\[
 \nu_{\boldsymbol\beta}
 =\delta_1+n_{\boldsymbol\beta}(u)\,du.
\]
There are constants $C_p>0$ and $C_{\boldsymbol\beta,p}>0$ such that
\begin{equation}\label{eq:image-density-growth}
 \boxed{
 |n_{\boldsymbol\beta}(u)|
 \le C_{\boldsymbol\beta,p}
 \exp\!\left(
 C_pM(\boldsymbol\beta)u\right),
 \qquad
 M(\boldsymbol\beta)
 =\max_{2\le j\le p}|\beta_j|^{2/(2j-1)}.}
\end{equation}
For $\boldsymbol\beta=\mathbf0$, one has
$\nu_{\mathbf0}=\delta_1$.
\end{theorem}

\begin{proof}
Choose $z_p>0$ with $\Phi_p(z_p)<1/4$.  If
$\boldsymbol\beta\ne\mathbf0$, choose $C_p$ sufficiently large that, with
$\sigma=C_pM(\boldsymbol\beta)$,
\[
 \sum_{j=2}^pC_p^{-(j-1/2)}\le z_p.
\]
Then $Z\le z_p$.  Choose $V\ge1$ so that
$(\sigma V)^{-1/2}\Psi_p(z_p)<1/4$.  By
\eqref{eq:defect-tail-norm}, the tail defect has norm smaller than $1/2$.

It remains to control the factorized forcing.  For every $\sigma_0>0$,
termwise use of
\[
 \frac{s^{q-1}}{\Gamma(q)}
 \le\sigma_0^{1-q}e^{\sigma_0s}
\]
in \eqref{eq:polynomial-image-kernel} gives
\begin{equation}\label{eq:kernel-direct-bound}
 |\mathcal K_{\boldsymbol\beta}(u,v)|
 \le\sigma_0e^{\sigma_0(u-v)}
 \left[
 \exp\!\left(
 \sqrt v\sum_{j=2}^p
 \frac{|\beta_j|}{2^{j-1}\sigma_0^{j-1}}
 \right)-1
 \right].
\end{equation}
Taking $0<\sigma_0<\sigma$ shows that both the single-kernel and composition
parts of $(\mathcal T^{\rm im}_{-\boldsymbol\beta})^2\delta_1$ belong to
$X_{\sigma,V}$.  Indeed, if
\[
 A_{\boldsymbol\beta}(\sigma_0)
 =\sum_{j=2}^p
 \frac{|\beta_j|}{2^{j-1}\sigma_0^{j-1}},
\]
the composition term, after multiplication by $e^{-\sigma u}$, is bounded
by a constant times
\[
 (1+u)\exp\!\left(
 -(\sigma-\sigma_0)u
 +A_{\boldsymbol\beta}(\sigma_0)\sqrt u\right),
\]
which is bounded on $[1,\infty)$.  Standard Volterra theory gives the
solution on $[1,V]$ \cite{JKV,Peskir};
moving the compact history into the forcing and inverting the tail by a
Neumann series extends it uniquely to $[1,\infty)$.  Finally, if
$R=\mathcal T^{\rm im}_{\boldsymbol\beta}\nu_{\boldsymbol\beta}
-\mathcal T^{\rm im}_{-\boldsymbol\beta}\delta_1$, the factorized equation
gives $\mathcal T^{\rm im}_{-\boldsymbol\beta}R=0$.  Uniqueness for the
homogeneous Volterra equation implies $R=0$ on every compact interval.
\end{proof}

\section{The killed heat kernel and first-passage density}
\label{sec:killed-kernel}

Let $\mu_{\bc,a}$ be the push-forward of
$\nu_{\boldsymbol\beta}$ under
\[
 u\longmapsto d+a\sqrt u.
\]
Then \eqref{eq:image-potential} becomes
\begin{align}
 P_{\bc,a}(t,x)
 ={}&\G(t,x-d+a)-\G(t,x-d-a)
 \notag\\
 &-\int_1^\infty
 \G(t,x-d-a\sqrt u)n_{\boldsymbol\beta}(u)\,du.
 \label{eq:polynomial-killed-potential}
\end{align}
Put
\[
 M(\bc)=\max_{2\le j\le p}|c_j|^{2/(2j-1)},
\]
and let $T_{\rm geom}$ be the first positive time, if any, at which
$F(t)=x_0$, with $T_{\rm geom}=\infty$ if no such time exists.

\begin{theorem}[Polynomial first passage]
\label{thm:polynomial-first-passage}
There exists $c_p>0$, depending only on $p$, such that
\eqref{eq:polynomial-killed-potential}, $\partial_xP_{\bc,a}$,
$\partial_{xx}P_{\bc,a}$, and $\partial_tP_{\bc,a}$ converge locally
uniformly whenever
\[
 0<t<T_{\rm an}:=\frac{c_p}{M(\bc)},
\]
with $T_{\rm an}=\infty$ when $\bc=\mathbf0$.  On
\begin{equation}\label{eq:polynomial-valid-time}
 0<t<T_{\bc,a,p}:=T_{\rm an}\wedge T_{\rm geom},
\end{equation}
$P_{\bc,a}$ is the transition density of Brownian motion started at $x_0$
and killed at $F$.  The first-passage density is
\begin{equation}\label{eq:polynomial-flux-density}
 q_{\bc,a}(t)=-\frac12
 \partial_xP_{\bc,a}(t,F(t)-),
\end{equation}
namely
\begin{align}
 q_{\bc,a}(t)
 ={}&\frac{a+f_{\bc}(t)}{2t}
 \G(t,a+f_{\bc}(t))
 -\frac{f_{\bc}(t)-a}{2t}
 \G(t,f_{\bc}(t)-a)
 \notag\\
 &-\frac1{2t}\int_1^\infty
 (f_{\bc}(t)-a\sqrt u)
 \G(t,f_{\bc}(t)-a\sqrt u)
 n_{\boldsymbol\beta}(u)\,du.
 \label{eq:polynomial-density-explicit}
\end{align}
\end{theorem}

\begin{proof}
Theorem~\ref{thm:polynomial-image-density} and
$M(\boldsymbol\beta)=a^2M(\bc)$ show that the Gaussian in
\eqref{eq:polynomial-killed-potential} dominates the image tail whenever
$a^2/(2t)>C_pM(\boldsymbol\beta)$.  Each displayed derivative introduces
only a polynomial factor in $\sqrt u$, so the same strict inequality gives
local uniform convergence.  Decreasing $c_p$ gives the stated interval.
At $x=F(t)$, the zero trace is exactly
\eqref{eq:polynomial-image-transform}.  The initial datum inside
$(-\infty,d)$ is $\delta_{x_0}$: every image point lies in
$[d+a,\infty)$ and hence contributes zero in the distributional limit on
compact subsets of $(-\infty,d)$.  On \eqref{eq:polynomial-valid-time}, the
boundary is smooth and remains strictly above $x_0$.  The Gaussian estimate
also gives decay as $x\to-\infty$.  Uniqueness for the moving Dirichlet heat
problem, by the maximum principle, identifies $P_{\bc,a}$ with the killed
Brownian kernel; see also \cite{Lee} for first-passage regularity under
substantially weaker boundary assumptions.  Differentiating its survival
probability gives \eqref{eq:polynomial-flux-density}, and direct
differentiation gives \eqref{eq:polynomial-density-explicit}.
\end{proof}

\section{Polynomial grading and the transported-zero construction}

For $\bn=(\nu_2,\ldots,\nu_p)\in\N_0^{p-1}$, write
\begin{equation}\label{eq:multi-notation}
 |\bn|=\sum_{j=2}^p\nu_j,
 \quad
 \bn!=\prod_{j=2}^p\nu_j!,
 \quad
 \bc^{\bn}=\prod_{j=2}^pc_j^{\nu_j},
\end{equation}
and define
\begin{equation}\label{eq:KW}
 K(\bn)=\sum_{j=2}^pj\nu_j,
 \qquad
 W(\bn)=\sum_{j=2}^p(2j-1)\nu_j.
\end{equation}
Then
\begin{equation}\label{eq:balance}
 W(\bn)=2K(\bn)-|\bn|.
\end{equation}

\begin{definition}\label{def:A}
For $\bn\in\N_0^{p-1}$ define $\A_{\bn}$ on monomials by
\begin{equation}\label{eq:A}
 \A_{\bn}x^n=C_{\bn}(n)x^{n+W(\bn)},
\end{equation}
where $\A_{\mathbf0}=I$ and
\begin{equation}\label{eq:C}
 C_{\bn}(n)=
 \frac{(-1)^{|\bn|}}{\bn!}
 \frac{\displaystyle\prod_{\ell=0}^{K(\bn)-1}(n+2\ell+1)}
 {\displaystyle\prod_{\ell=1}^{W(\bn)}(n+\ell)}.
\end{equation}
The polynomial spatial transmutation is the formal series
\begin{equation}\label{eq:T}
 \T_{\bc}=\sum_{\bn\ge\mathbf0}\bc^{\bn}\A_{\bn}.
\end{equation}
\end{definition}

The following theorem gives the transported-zero origin of
Definition~\ref{def:A}.

\begin{theorem}[Polynomial transport of odd heat modes]
\label{thm:odd-transport}
Let $n=2h+1$ be odd and define
\[
 \Gamma_{n,\bc}=\T_{\bc}x^n.
\]
Then, as formal series in $\bc$ and $t$,
\begin{align}
 (P_t\Gamma_{n,\bc})(F_{\bc}(t))&=0,
 \label{eq:odd-zero}\\
 \partial_x(P_t\Gamma_{n,\bc})(F_{\bc}(t))
 &=\partial_x(P_tx^n)(0)
 =\frac{n!}{2^hh!}t^h.
 \label{eq:odd-flux}
\end{align}
\end{theorem}

\begin{proof}
For $d\ge0$,
\begin{equation}\label{eq:heat-monomial}
 P_tx^d=\sum_{r=0}^{\lfloor d/2\rfloor}
 \frac{d!}{(d-2r)!r!}\left(\frac t2\right)^rx^{d-2r}.
\end{equation}
Fix a target multi-index $\bm$.  A contribution from $\bn\le\bm$ to the
coefficient of $\bc^{\bm}$ in the boundary trace must use
$|\bm|-|\bn|$ factors from $F_{\bc}$.  Consequently the heat index is
\[
 r=h-\frac{|\bm|-1}{2}+K(\bn).
\]
If $|\bm|$ is even, the required heat index is not integral and the
coefficient vanishes.  If $|\bm|=2L+1$, substitution of \eqref{eq:C} and
\eqref{eq:heat-monomial} reduces the coefficient, up to a factor independent
of $\bn$, to
\begin{equation}\label{eq:finite-difference}
 \sum_{\bn\le\bm}(-1)^{|\bn|}
 \binom{\bm}{\bn}
 (K(\bn)+h)_{\underline L},
\end{equation}
where $(z)_{\underline L}=z(z-1)\cdots(z-L+1)$.  The expression in
\eqref{eq:finite-difference} is the full multivariate finite difference of a
polynomial of total degree $L$.  Since $L<|\bm|$, it is zero.
Terms for which the displayed heat index is negative are harmless in this
reduction: then $K(\bn)+h<L$ and the falling factorial vanishes.

For the spatial derivative, a nonzero coefficient requires $|\bm|=2L$.
The same cancellation reduces it to the full finite difference of
$(K(\bn)+h)_{\underline L}$, now with $L<|\bm|$ whenever $\bm\ne\mathbf0$.
Thus every nonconstant coefficient vanishes.  At $\bm=\mathbf0$ one obtains
the flat flux
\[
 \partial_x(P_tx^{2h+1})(0)=\frac{(2h+1)!}{2^hh!}t^h.
\]
This proves both identities.
\end{proof}

Taking $n=1$ gives the distinguished constant-flux seed.

\begin{corollary}[Entire constant-flux selection]
\label{cor:Gamma}
The function
\begin{equation}\label{eq:Gamma}
 \boxed{
 \Gamma_{\bc}(x)=\T_{\bc}x
 =\sum_{\bn\ge\mathbf0}
 (-1)^{|\bn|}2^{K(\bn)}
 \frac{K(\bn)!}{\bn!\,[W(\bn)+1]!}
 \bc^{\bn}x^{W(\bn)+1}}
\end{equation}
is entire in $x$ for every fixed finite vector $\bc$.  Formally,
\begin{equation}\label{eq:Gamma-identities}
 (P_t\Gamma_{\bc})(F_{\bc}(t))=0,
 \qquad
 \partial_x(P_t\Gamma_{\bc})(F_{\bc}(t))=1.
\end{equation}
\end{corollary}

\begin{proof}
The coefficient formula follows from
\[
 \prod_{\ell=0}^{K-1}(2\ell+2)=2^KK!,
 \qquad
 \prod_{\ell=1}^{W}(1+\ell)=(W+1)!.
\]
To see entireness, group the terms by $R=|\bn|$.  Since
$2R\le K(\bn)\le pR$ and $W(\bn)=2K(\bn)-R$, one has
\[
 \frac{K(\bn)!}{[W(\bn)+1]!}
 =\frac{K(\bn)!}{[2K(\bn)-R+1]!}
 \le \frac1{(2R+1)^{R+1}}.
\]
Moreover,
\[
 \sum_{|\bn|=R}\frac1{\bn!}=\frac{(p-1)^R}{R!}.
\]
On every compact $x$-set, all other factors in \eqref{eq:Gamma} are bounded
by a constant to the power $R$.  Hence the level-$R$ contribution is bounded
by
\[
 \frac{C^R}{R!(2R+1)^{R+1}}
\]
for a compact-dependent constant $C$.  The grouped series therefore
converges normally on every compact subset of $\mathbb C$.
\end{proof}

\begin{remark}
Constant and linear terms can be restored separately.  A spatial translation
handles a constant term, while the Galilean symmetry of the heat equation
handles a linear term.  We work with \eqref{eq:F} in order to isolate the
nonlinear polynomial structure.
\end{remark}

\section{Laplace conjugacy}

Let
\begin{equation}\label{eq:Laplace}
 (\Lap f)(k)=\int_0^\infty e^{-kx}f(x)\,dx,
 \qquad
 \Dt=-\frac1k\frac{d}{dk}.
\end{equation}
Algebraically, $\Lap(x^n)=n!k^{-n-1}$.

\begin{theorem}[Multi-index Laplace conjugacy]
\label{thm:conjugacy}
For every $\bn\in\N_0^{p-1}$,
\begin{equation}\label{eq:conjugacy}
 \boxed{
 \Lap\A_{\bn}\Lap^{-1}
 =\frac{(-1)^{|\bn|}}{\bn!}
 k^{|\bn|}\Dt^{K(\bn)}.}
\end{equation}
The identity holds algebraically on monomials and formal power series, and on
any analytic class for which the transformations are justified.
\end{theorem}

\begin{proof}
Since
\[
 \Dt^K(k^{-n-1})=
 \prod_{\ell=0}^{K-1}(n+2\ell+1)k^{-n-1-2K},
\]
the right-hand side of \eqref{eq:conjugacy}, applied to $\Lap(x^n)$, is
\[
 \frac{(-1)^{|\bn|}}{\bn!}n!
 \prod_{\ell=0}^{K(\bn)-1}(n+2\ell+1)
 k^{-n-1-2K(\bn)+|\bn|}.
\]
By \eqref{eq:balance}, the exponent is $-n-W(\bn)-1$.  On the other hand,
\[
 \Lap(\A_{\bn}x^n)
 =C_{\bn}(n)(n+W(\bn))!k^{-n-W(\bn)-1}.
\]
Using
\[
 \frac{(n+W)!}{\prod_{\ell=1}^W(n+\ell)}=n!
\]
gives the same expression.
\end{proof}

\begin{corollary}\label{cor:sign-reversed}
The sign-reversed transmutation satisfies
\begin{equation}\label{eq:T-conjugacy}
 \boxed{
 \Lap\T_{-\bc}\Lap^{-1}
 =\sum_{\bn\ge\mathbf0}
 \frac{\bc^{\bn}}{\bn!}
 k^{|\bn|}\Dt^{K(\bn)}.}
\end{equation}
\end{corollary}

\begin{remark}[Ordering and interaction]
The operator in \eqref{eq:T-conjugacy} is ordered: multiplication by
$k^{|\bn|}$ occurs after $\Dt^{K(\bn)}$.  It is not an ordinary exponential
of a sum of commuting monomial generators.  Correspondingly, mixed spatial
operators are generally not compositions of pure ones.  For distinct
$r,s\in\{2,\ldots,p\}$,
\[
 \A_{\mathbf e_r+\mathbf e_s}
 \ne \A_{\mathbf e_r}\A_{\mathbf e_s},
 \qquad
 \A_{\mathbf e_r+\mathbf e_s}
 \ne \A_{\mathbf e_s}\A_{\mathbf e_r}.
\]
The polynomial theory therefore contains genuine interactions.
\end{remark}

\section{The formal first-passage theorem}

Let
\[
 \G(t,y)=\frac1{\sqrt{2\pi t}}e^{-y^2/(2t)}.
\]
For a function or distribution $g$ supported in $[a,\infty)$, define the
flat flux map, with the integral interpreted as a distributional pairing,
\begin{equation}\label{eq:B}
 (\Bf g)(t)=\frac1t\left\langle g(y),y\G(t,y)\right\rangle.
\end{equation}
The classical identity
\[
 \int_0^\infty e^{-k^2t/2}\frac{y}{t}\G(t,y)\,dt=e^{-ky}
\]
gives
\begin{equation}\label{eq:flux-laplace}
 \int_0^\infty e^{-k^2t/2}(\Bf g)(t)\,dt=(\Lap g)(k).
\end{equation}

For $k>0$, the process
\[
 \exp\!\left(kB_t-\frac{k^2t}{2}\right)
\]
is a martingale.  Under the usual integrability and terminal conditions,
optional stopping at \eqref{eq:tau} yields
\begin{equation}\label{eq:martingale-transform}
 \int_0^\infty
 e^{-k^2t/2+kF_{\bc}(t)}q_{\bc,a}(t)\,dt=e^{-ak},
\end{equation}
where $q_{\bc,a}$ denotes the first-passage density.  We use this identity
formally at $\bc=\mathbf0$.  Theorem~\ref{thm:polynomial-first-passage}
constructs the genuine density on a nonzero interval, but convergence of its
global martingale transform and termwise differentiation in $\bc$ are
separate issues.

Introduce
\begin{equation}\label{eq:qformal}
 q_{\bc,a}^{\mathrm{form}}(t)
 =\sum_{\bm\ge\mathbf0}\bc^{\bm}q_{\bm}(t;a),
 \qquad
 \widehat q_{\bm}(k)=
 \int_0^\infty e^{-k^2t/2}q_{\bm}(t;a)\,dt.
\end{equation}
Because
\begin{equation}\label{eq:expansion-exponential}
 e^{kF_{\bc}(t)}
 =\sum_{\bn\ge\mathbf0}
 \frac{k^{|\bn|}}{\bn!}\bc^{\bn}t^{K(\bn)}
\end{equation}
and
\[
 \Dt^K\widehat q(k)
 =\int_0^\infty t^Ke^{-k^2t/2}q(t)\,dt,
\]
coefficient extraction gives the following recursion.

\begin{proposition}[Martingale coefficient recursion]
\label{prop:mart-recursion}
The formal coefficients in \eqref{eq:qformal} are uniquely determined by
\begin{equation}\label{eq:mart-recursion}
 \boxed{
 \sum_{\bn\le\bm}
 \frac{k^{|\bn|}}{\bn!}\Dt^{K(\bn)}
 \widehat q_{\bm-\bn}(k)
 =\mathbf1_{\{\bm=\mathbf0\}}e^{-ak}.}
\end{equation}
\end{proposition}

We now obtain the same recursion from the spatial theory.  Formally invert
\begin{equation}\label{eq:inverse-T}
 \T_{-\bc}^{-1}
 =\sum_{\bm\ge\mathbf0}\bc^{\bm}\mathcal S_{\bm}.
\end{equation}
Here the action of $\A_{\bn}$ on a distribution supported in
$[a,\infty)$ is defined operationally by
\begin{equation}\label{eq:A-distributional}
 \Lap(\A_{\bn}g)
 =\frac{(-1)^{|\bn|}}{\bn!}
 k^{|\bn|}\Dt^{K(\bn)}\Lap g.
\end{equation}
For $g=\delta_a$, and recursively for all distributions used below, the
right-hand side is the Laplace transform of a distribution supported in
$[a,\infty)$.  Thus \eqref{eq:A-distributional} extends the algebraic
definition without assuming that $\A_{\bn}$ is a finite-order differential
operator.
Then
\begin{equation}\label{eq:S-recursion}
 \mathcal S_{\mathbf0}=I,
 \qquad
 \mathcal S_{\bm}
 =\sum_{\mathbf0<\bn\le\bm}
 (-1)^{|\bn|+1}\A_{\bn}\mathcal S_{\bm-\bn}.
\end{equation}
Set
\begin{equation}\label{eq:gh}
 g_{\bm}=\mathcal S_{\bm}\delta_a,
 \qquad
 h_{\bm}=\Bf g_{\bm}.
\end{equation}

\begin{theorem}[Formal polynomial first-passage transmutation]
\label{thm:first-passage}
For every finite polynomial boundary \eqref{eq:F}, the flux coefficients
generated by \eqref{eq:inverse-T}--\eqref{eq:gh} coincide with the unique
formal martingale sequence:
\begin{equation}\label{eq:hq}
 \boxed{h_{\bm}(t;a)=q_{\bm}(t;a),
 \qquad \bm\in\N_0^{p-1}.}
\end{equation}
Equivalently,
\begin{equation}\label{eq:h-recursion}
 \sum_{\bn\le\bm}
 \frac{k^{|\bn|}}{\bn!}\Dt^{K(\bn)}
 \widehat h_{\bm-\bn}(k)
 =\mathbf1_{\{\bm=\mathbf0\}}e^{-ak}.
\end{equation}
\end{theorem}

\begin{proof}
The inverse identity is
\[
 \T_{-\bc}\left(
 \sum_{\bm\ge\mathbf0}\bc^{\bm}g_{\bm}
 \right)=\delta_a.
\]
Apply $\Lap$ and use Corollary~\ref{cor:sign-reversed} to obtain
\[
 \left(
 \sum_{\bn\ge\mathbf0}
 \frac{\bc^{\bn}}{\bn!}k^{|\bn|}\Dt^{K(\bn)}
 \right)
 \left(
 \sum_{\bm\ge\mathbf0}\bc^{\bm}\Lap g_{\bm}
 \right)=e^{-ak}.
\]
By \eqref{eq:flux-laplace}, $\Lap g_{\bm}=\widehat h_{\bm}$.  Equating
coefficients gives \eqref{eq:h-recursion}, which is identical to
\eqref{eq:mart-recursion}.  The recursion is triangular with leading term
$\widehat h_{\bm}$, so its solution is unique.
\end{proof}

\begin{remark}
Theorem~\ref{thm:first-passage} is formal at all orders.  It asserts neither
convergence of $\sum\bc^{\bm}h_{\bm}$ nor analytic dependence of the genuine
density on $\bc$.  Whenever such an expansion of the genuine density can be
justified, however, its coefficients must equal those generated by the
inverse transmutation.
\end{remark}

\section{The quadratic--cubic boundary}
\label{sec:quadratic-cubic}

Consider the first genuinely mixed case
\begin{equation}\label{eq:AB-boundary}
 F_{A,B}(t)=At^2+Bt^3.
\end{equation}
For $\bn=(m,n)$,
\[
 K(m,n)=2m+3n,
 \qquad
 W(m,n)=3m+5n.
\]
The constant-flux function is
\begin{equation}\label{eq:GammaAB}
 \boxed{
 \Gamma_{A,B}(x)=
 x\sum_{m,n\ge0}
 \frac{(2m+3n)!}{m!n!(3m+5n+1)!}
 (-4Ax^3)^m(-8Bx^5)^n.}
\end{equation}
Its first terms are
\begin{align}\label{eq:GammaAB-first}
 \Gamma_{A,B}(x)
 ={}&x-\frac A3x^4-\frac B{15}x^6
 +\frac{4A^2}{105}x^7+\frac{2AB}{189}x^9
 +\frac{2B^2}{3465}x^{11}+\cdots.
\end{align}
Theorem~\ref{thm:odd-transport} gives
\[
 (P_t\Gamma_{A,B})(At^2+Bt^3)=0,
 \qquad
 \partial_x(P_t\Gamma_{A,B})(At^2+Bt^3)=1.
\]

The spatial operators are
\begin{equation}\label{eq:Amn}
 \A_{m,n}x^r=
 \frac{(-1)^{m+n}}{m!n!}
 \frac{\prod_{\ell=0}^{2m+3n-1}(r+2\ell+1)}
 {\prod_{\ell=1}^{3m+5n}(r+\ell)}
 x^{r+3m+5n},
\end{equation}
and
\begin{equation}\label{eq:Amn-conjugacy}
 \Lap\A_{m,n}\Lap^{-1}
 =\frac{(-1)^{m+n}}{m!n!}k^{m+n}\Dt^{2m+3n}.
\end{equation}

The double series \eqref{eq:GammaAB} is of Horn type.  Put
\[
 H(z,w)=\frac{\Gamma_{A,B}(x)}x,
 \qquad z=-4Ax^3,
 \qquad w=-8Bx^5,
\]
and let $\theta_z=z\partial_z$, $\theta_w=w\partial_w$,
$Q=3\theta_z+5\theta_w$, $K=2\theta_z+3\theta_w$.  Direct coefficient
comparison gives the explicit Horn system
\begin{align}
 \theta_z(Q-1)Q(Q+1)H
 &=z(K+1)(K+2)H,
 \label{eq:Horn-z}\\
 \theta_w(Q-3)(Q-2)(Q-1)Q(Q+1)H
 &=w(K+1)(K+2)(K+3)H.
 \label{eq:Horn-w}
\end{align}
Thus the mixed constant-flux function is holonomic.  The elimination can in
fact be carried out explicitly.  Write
\[
 (\theta)_{\underline r}=\theta(\theta-1)\cdots(\theta-r+1),
 \qquad \theta=xD.
\]

\begin{theorem}[A scalar quadratic--cubic ODE]
\label{thm:mixed-ode}
The function $\Gamma_{A,B}$ in \eqref{eq:GammaAB} is annihilated by the
fifth-order Euler-form operator
\begin{equation}\label{eq:mixed-operator}
 \mathcal L_{A,B}=\sum_{(i,j)\in\mathcal S}
 A^iB^jx^{3i+5j}P_{ij}(\theta),
\end{equation}
where
\[
 \mathcal S=\{(0,4),(0,5),(1,4),(2,3),(2,4),(3,3),
 (4,2),(5,1),(5,2),(6,1),(7,0),(8,0)\}
\]
and
\begin{align*}
 P_{04}(\theta)&=16875(\theta)_{\underline5},\\
 P_{05}(\theta)&=16200(3\theta+2)(3\theta+7),\\
 P_{14}(\theta)&=48600(\theta)_{\underline2}(2\theta+1),\\
 P_{23}(\theta)&=675(\theta)_{\underline4}(5\theta-31),\\
 P_{24}(\theta)&=3240(3\theta+2)(3\theta+4),\\
 P_{33}(\theta)&=1080\theta(18\theta^2-47\theta-26),\\
 P_{42}(\theta)&=-90(\theta)_{\underline3}(9\theta-26),\\
 P_{51}(\theta)&=50(\theta)_{\underline5},\\
 P_{52}(\theta)&=-48(81\theta^2+27\theta-8),\\
 P_{61}(\theta)&=24(\theta)_{\underline2}(12\theta+5),\\
 P_{70}(\theta)&=-12(\theta)_{\underline4},\\
 P_{80}(\theta)&=-32\theta(2\theta+1).
\end{align*}
Thus
\begin{equation}\label{eq:mixed-ode}
 \boxed{\mathcal L_{A,B}\Gamma_{A,B}=0.}
\end{equation}
No assertion of minimality is made.
\end{theorem}

\begin{proof}
Let
\[
 c_{m,n}=
 \frac{(-1)^{m+n}2^{2m+3n}(2m+3n)!}
 {m!n!(3m+5n+1)!},
\]
so that
\[
 \Gamma_{A,B}(x)=\sum_{m,n\ge0}
 c_{m,n}A^mB^nx^{3m+5n+1}.
\]
Put $K=2m+3n$ and $W=3m+5n$.  For fixed nonnegative $i,j$,
\begin{equation}\label{eq:coefficient-ratio}
 \frac{c_{m-i,n-j}}{c_{m,n}}
 =\frac{(-1)^{i+j}}{2^{2i+3j}}
 (m)_{\underline i}(n)_{\underline j}
 \frac{(W+1)_{\underline{3i+5j}}}
 {(K)_{\underline{2i+3j}}},
\end{equation}
with the convention that a term vanishes when $i>m$ or $j>n$.
The coefficient of $A^mB^nx^{W+1}$ in
$\mathcal L_{A,B}\Gamma_{A,B}$, divided by $c_{m,n}$, is therefore
\begin{equation}\label{eq:ode-coefficient-identity}
 \sum_{\substack{(i,j)\in\mathcal S\\i\le m,\,j\le n}}
 \frac{c_{m-i,n-j}}{c_{m,n}}
 P_{ij}(W+1-3i-5j).
\end{equation}
Substitution of \eqref{eq:coefficient-ratio} and the twelve displayed
polynomials, followed by clearing the falling-factorial denominators,
reduces \eqref{eq:ode-coefficient-identity} to the zero polynomial in
$m,n$.  Hence every coefficient vanishes.  This is an exact rational
identity; no numerical interpolation is used.
\end{proof}

The coordinate-axis degenerations contain the monomial equations in a
particularly transparent way.  Define
\begin{align}
 \mathcal M_A
 &=3(\theta)_{\underline3}+8Ax^3(2\theta+1),
 \label{eq:quadratic-minimal}\\
 \mathcal M_B
 &=25(\theta)_{\underline5}
 +24Bx^5(3\theta+2)(3\theta+7).
 \label{eq:cubic-minimal}
\end{align}
These are the distinguished order-three quadratic and order-five cubic
operators.  Direct
use of
\[
 P(\theta)x^r=x^rP(\theta+r)
\]
gives the factorizations
\begin{align}
 \mathcal L_{A,0}
 ={}&-4A^7x^{21}(\theta-3)\circ\mathcal M_A,
 \label{eq:A-axis-factorization}\\
 \mathcal L_{0,B}
 ={}&675B^4x^{20}\mathcal M_B.
 \label{eq:B-axis-factorization}
\end{align}
Therefore the mixed equation specializes to left multiples of the correct
monomial equations.  The extra factors record the degeneration of the mixed
system on the coordinate axes; they do not alter the orders of the displayed
monomial annihilators.

\subsection{A compact low-valuation annihilator}
\label{subsec:order-eight}

Lower differential order is not the only useful notion of simplicity.  The
same function is annihilated by the order-eight operator
\begin{equation}\label{eq:L8}
 \widetilde{\mathcal L}_{A,B}
 =\sum_{(i,j)\in\widetilde{\mathcal S}}
 A^iB^jx^{3i+5j}\widetilde P_{ij}(\theta),
\end{equation}
where
\[
 \widetilde{\mathcal S}
 =\{(0,1),(0,2),(1,1),(1,2),(2,0),(2,1),(3,0),(4,0)\}
\]
and
\begin{align*}
 \widetilde P_{01}&=-125(\theta)_{\underline8},\\
 \widetilde P_{02}&=-120(\theta)_{\underline3}(3\theta+2)(3\theta+7),\\
 \widetilde P_{11}&=-40(\theta)_{\underline5}(33\theta-46),\\
 \widetilde P_{12}&=-192(3\theta+2)(3\theta+4)(3\theta+7),\\
 \widetilde P_{20}&=30(\theta)_{\underline7},\\
 \widetilde P_{21}&=-96(\theta)_{\underline2}(36\theta^2-3\theta-10),\\
 \widetilde P_{30}&=16(\theta)_{\underline4}(19\theta-10),\\
 \widetilde P_{40}&=128\theta(2\theta+1)(3\theta+4).
\end{align*}
The same coefficient-ratio calculation used in
Theorem~\ref{thm:mixed-ode} proves
\begin{equation}\label{eq:L8-annihilates}
 \widetilde{\mathcal L}_{A,B}\Gamma_{A,B}=0.
\end{equation}
Its coordinate-axis restrictions are
\begin{align}
 \widetilde{\mathcal L}_{A,0}
 ={}&A^2\Bigl[
 10x^6(\theta-3)(\theta-4)(\theta-5)(\theta-6)
 \notag\\[-1mm]
 &\hspace{29mm}
 +16Ax^9(\theta-3)(3\theta-5)
 \Bigr]\circ\mathcal M_A,
 \label{eq:L8-A-axis}\\
 \widetilde{\mathcal L}_{0,B}
 ={}&-5Bx^5(\theta-5)(\theta-6)(\theta-7)
 \circ\mathcal M_B.
 \label{eq:L8-B-axis}
\end{align}

For an Euler-form operator
$L=\sum A^iB^jx^{3i+5j}P_{ij}(\theta)$, define its parameter valuation by
\[
 \operatorname{val}_{A,B}(L)
 =\min\{i+j:P_{ij}\ne0\}.
\]
The two annihilators have the complementary complexity profiles
\begin{center}
\begin{tabular}{lccc}
\toprule
operator & differential order & parameter valuation & stencil size\\
\midrule
$\mathcal L_{A,B}$ & $5$ & $4$ & $12$\\
$\widetilde{\mathcal L}_{A,B}$ & $8$ & $1$ & $8$\\
\bottomrule
\end{tabular}
\end{center}
Thus $\mathcal L_{A,B}$ has lower differential order, whereas
$\widetilde{\mathcal L}_{A,B}$ is visible at much lower order in a formal
expansion about $(A,B)=(0,0)$.

The distinction is algebraic rather than accidental.  Over the rational Ore
field $\mathbb Q(A,B,x)\langle D\rangle$, exact right division gives
\begin{equation}\label{eq:L8-factor-L5}
 \boxed{
 \widetilde{\mathcal L}_{A,B}
 =\mathcal Q_{A,B}\circ\mathcal L_{A,B},}
\end{equation}
where $\mathcal Q_{A,B}$ has order three.  Its coefficients are recorded in
Appendix~\ref{app:quotient}.  In particular, the leading coefficient is
\[
 [D^3]\mathcal Q_{A,B}
 =-\frac5{x^{12}\Delta(A,B,x)},
 \qquad
 \Delta(A,B,x)=2A^5+135A^2B^2x+675B^3.
\]
Consequently, over the generic rational field the order-eight annihilation
condition carries no independent constraint on the distinguished solution
beyond the order-five equation.  This is not an equivalence of solution
spaces: the order-eight operator has a larger formal kernel.  Moreover, the
quotient is singular at the perturbative base point $A=B=0$, so the two
operators need not be equally useful over $\mathbb Q[A,B,x]$ or over formal
parameter rings.

This observation is relevant to first-passage applications.  The density
construction in Theorem~\ref{thm:first-passage} comes from the full spatial
transmutation, not from the order of a single annihilating ODE.  Moreover,
neither of the two displayed operators annihilates all transported odd modes;
for example, direct substitution already gives nonzero coefficients when
they act on $\T_{A,B}x^3$.  A higher-order annihilator therefore supplies a
larger formal kernel, but its additional elements are not automatically
adapted to the moving-zero boundary condition.

\subsection{The coefficient matrix and the polynomial hierarchy}

The preceding calculation has a direct matrix formulation that extends to
an arbitrary finite polynomial boundary.  Write
\begin{equation}\label{eq:general-coefficients}
 \Gamma_{\bc}(x)=\sum_{\bm\ge\mathbf0}
 \gamma_{\bm}\bc^{\bm}x^{W(\bm)+1},
 \qquad
 \gamma_{\bm}=(-1)^{|\bm|}2^{K(\bm)}
 \frac{K(\bm)!}{\bm!\,[W(\bm)+1]!}.
\end{equation}
Fix a finite multi-index stencil $\mathcal S$ and an integer $d\ge0$, and
seek an Euler-form operator
\begin{equation}\label{eq:general-euler-ansatz}
 \mathcal L=\sum_{\bn\in\mathcal S}
 \bc^{\bn}x^{W(\bn)}P_{\bn}(\theta),
 \qquad
 P_{\bn}(z)=\sum_{r=0}^dp_{\bn,r}z^r.
\end{equation}
Define the coefficient matrix $\mathsf M_{\mathcal S,d}$, with rows indexed
by $\bm\in\N_0^{p-1}$ and columns indexed by $(\bn,r)$, by
\begin{equation}\label{eq:coefficient-matrix}
 \boxed{
 (\mathsf M_{\mathcal S,d})_{\bm,(\bn,r)}
 =\mathbf1_{\{\bn\le\bm\}}
 \gamma_{\bm-\bn}
 [W(\bm-\bn)+1]^r.}
\end{equation}

\begin{proposition}[Matrix characterization of polynomial ODEs]
\label{prop:matrix-odes}
The operator \eqref{eq:general-euler-ansatz} annihilates $\Gamma_{\bc}$ if
and only if its coefficient vector
$\boldsymbol p=(p_{\bn,r})$ satisfies
\begin{equation}\label{eq:matrix-kernel}
 \boxed{\mathsf M_{\mathcal S,d}\boldsymbol p=0.}
\end{equation}
\end{proposition}

\begin{proof}
The coefficient of
$\bc^{\bm}x^{W(\bm)+1}$ in $\mathcal L\Gamma_{\bc}$ is
\[
 \sum_{\substack{\bn\in\mathcal S\\\bn\le\bm}}
 \gamma_{\bm-\bn}
 P_{\bn}(W(\bm-\bn)+1),
\]
which is exactly the $\bm$th component of
$\mathsf M_{\mathcal S,d}\boldsymbol p$.
\end{proof}

Increasing the stencil or the allowed polynomial degree produces nested
linear systems and, consequently, a hierarchy of annihilating operators.
Theorem~\ref{thm:mixed-ode} is obtained from
\[
 \mathcal S=
 \{(0,4),(0,5),(1,4),(2,3),(2,4),(3,3),
 (4,2),(5,1),(5,2),(6,1),(7,0),(8,0)\},
 \qquad d=5.
\]
Finite row truncations of \eqref{eq:coefficient-matrix} can be used to
discover candidate kernel vectors; the resulting identities must then be
verified for arbitrary $\bm$, as in the proof of
Theorem~\ref{thm:mixed-ode}.  This separates the constructive linear-algebra
step from the exact proof.

The matrix formulation also explains the coordinate-axis limits.  Setting
all parameters except $c_j$ equal to zero removes every column whose
multi-index is not supported on the $j$th axis.  Any nonzero surviving
operator annihilates the monomial constant-flux function.  In the present
quadratic--cubic case, exact Ore division gives the explicit left
factorizations \eqref{eq:A-axis-factorization}--\eqref{eq:B-axis-factorization}.
No general minimality assertion is needed for this conclusion.

The monomial equations and Theorem~\ref{thm:mixed-ode} suggest a broader
order principle.

\begin{remark}[Conjectural highest-degree hierarchy]
Let $c_p\ne0$.  We conjecture that the constant-flux function associated
with
\[
 F(t)=\sum_{j=2}^pc_jt^j
\]
admits a nonzero Euler-form annihilator of order at most $2p-1$.  On each
coordinate axis, a suitable member of the hierarchy should specialize to a
left multiple of the order-$2j-1$ monomial operator.  The monomial theory
proves the assertion on every coordinate axis, while
Theorem~\ref{thm:mixed-ode} proves the first genuinely mixed case
$(c_2,c_3)=(A,B)$ with order five.  Neither minimality nor the general
existence assertion is claimed here.
\end{remark}

A second useful representation is obtained by taking the spatial Laplace
transform:
\begin{equation}\label{eq:GammaAB-Laplace}
 \boxed{
 (\Lap\Gamma_{A,B})(q)
 =\frac1{q^2}\int_0^\infty
 \exp\!\left(-s-\frac{4As^2}{q^3}
 -\frac{8Bs^3}{q^5}\right)\,ds.}
\end{equation}
Equation \eqref{eq:GammaAB-Laplace} gives an alternative finite moment system
for the same holonomic function.  When $AB\ne0$, the equation contains eight
genuinely mixed parameter blocks in addition to the four coordinate-axis
blocks and does not split into independent quadratic and cubic pieces.

The first-passage recursion becomes
\begin{equation}\label{eq:AB-recursion}
 \boxed{
 \sum_{i=0}^{m}\sum_{j=0}^{n}
 \frac{k^{i+j}}{i!j!}\Dt^{2i+3j}
 \widehat q_{m-i,n-j}(k)
 =\delta_{m0}\delta_{n0}e^{-ak}.}
\end{equation}
Through total degree two,
\begin{align}
 \widehat q_{0,0}(k)&=e^{-ak},
 \label{eq:q00}\\
 \widehat q_{1,0}(k)
 &=-e^{-ak}\left(\frac{a^2}{k}+\frac a{k^2}\right),
 \label{eq:q10}\\
 \widehat q_{0,1}(k)
 &=-e^{-ak}\left(\frac{a^3}{k^2}
 +\frac{3a^2}{k^3}+\frac{3a}{k^4}\right),
 \label{eq:q01}\\
 \widehat q_{2,0}(k)
 &=\frac{ae^{-ak}}{2k^5}
 \left(a^3k^3+2a^2k^2+ak+1\right),
 \label{eq:q20}\\
 \widehat q_{1,1}(k)
 &=\frac{ae^{-ak}}{k^7}
 \left(a^4k^4+5a^3k^3+11a^2k^2+15ak+15\right),
 \label{eq:q11}\\
 \widehat q_{0,2}(k)
 &=\frac{ae^{-ak}}{2k^9}
 \left(a^5k^5+9a^4k^4+39a^3k^3+108a^2k^2
 +207ak+207\right).
 \label{eq:q02}
\end{align}
In particular, $q_{1,1}$ is the first genuine quadratic--cubic interaction.

\section{Numerical illustration}
\label{sec:numerical}

The calculation in this section is based on the multi-index transmutation,
rather than on a choice of scalar annihilating ODE.  It therefore tests the
structure common to both the order-five and order-eight differential
representations: the mixed coefficient generated by the simultaneous
quadratic and cubic boundary terms.

We illustrate the formal expansion with
\begin{equation}\label{eq:numerical-parameters}
 a=1,
 \qquad A=0.05,
 \qquad B=0.01,
 \qquad 0<t\le2.
\end{equation}
The benchmark is obtained independently from Durbin's second-kind equation.
Writing $b(t)=a+F_{A,B}(t)$ and
$\varphi(z)=(2\pi)^{-1/2}e^{-z^2/2}$, the density satisfies
\begin{align}\label{eq:Durbin-equation}
 q(t)={}&
 \left(\frac{b(t)}t-b'(t)\right)
 \frac{\varphi(b(t)/\sqrt t)}{\sqrt t}\\
 &+\int_0^tq(s)
 \left[b'(t)-\frac{b(t)-b(s)}{t-s}\right]
 \frac{\varphi((b(t)-b(s))/\sqrt{t-s})}{\sqrt{t-s}}\,ds.
 \notag
\end{align}
We discretize \eqref{eq:Durbin-equation} on uniform grids.  The largest
difference at the reported times between $2000$ and $4000$ time steps is
$4.42\times10^{-8}$.

The formal approximations are
\begin{align*}
 q^{[0]}&=q_{0,0},\\
 q^{[1]}&=q_{0,0}+Aq_{1,0}+Bq_{0,1},\\
 q^{[2]}&=q^{[1]}+A^2q_{2,0}+ABq_{1,1}+B^2q_{0,2}.
\end{align*}
The coefficient functions are recovered by numerical inversion of
\eqref{eq:q00}--\eqref{eq:q02}, with Laplace variable $s=k^2/2$.

\begin{table}[ht]
\centering
\caption{First-passage density for $F(t)=0.05t^2+0.01t^3$ and $a=1$.}
\label{tab:numerical}
\begin{tabular}{cccccc}
\toprule
$t$ & image flux & Durbin benchmark & $q^{[0]}$ & $q^{[1]}$ & $q^{[2]}$\\
\midrule
$0.2$ & $0.36165795$ & $0.36165838$ & $0.36612456$ & $0.36163531$ & $0.36165846$\\
$0.5$ & $0.39926278$ & $0.39926287$ & $0.41510750$ & $0.39904405$ & $0.39926340$\\
$1.0$ & $0.21751849$ & $0.21751846$ & $0.24197072$ & $0.21672331$ & $0.21751039$\\
$1.5$ & $0.12644009$ & $0.12644005$ & $0.15559955$ & $0.12469831$ & $0.12638022$\\
$2.0$ & $0.07750707$ & $0.07750706$ & $0.10984782$ & $0.07429092$ & $0.07730179$\\
\bottomrule
\end{tabular}
\end{table}

The image column is obtained independently by solving
\eqref{eq:polynomial-volterra-density} on $1\le u\le64$ with 3600
trapezoidal cells and total Wright-series degree ten.  Its largest difference
from the Durbin benchmark in the table is $4.34\times10^{-7}$; the boundary
trace residual is below $2.50\times10^{-7}$.

At $t=2$, the total-degree-one absolute error is approximately
$3.22\times10^{-3}$, whereas the total-degree-two error is
$2.05\times10^{-4}$.  The interaction contribution is
\[
 ABq_{1,1}(2)=1.475\times10^{-3}.
\]
If this mixed term is omitted while the two pure degree-two terms are kept,
the error is approximately $1.68\times10^{-3}$.  Thus the mixed coefficient
improves the degree-two approximation by a factor of about eight at the end
of the interval.

\begin{figure}[ht]
\centering
\includegraphics[width=\textwidth]{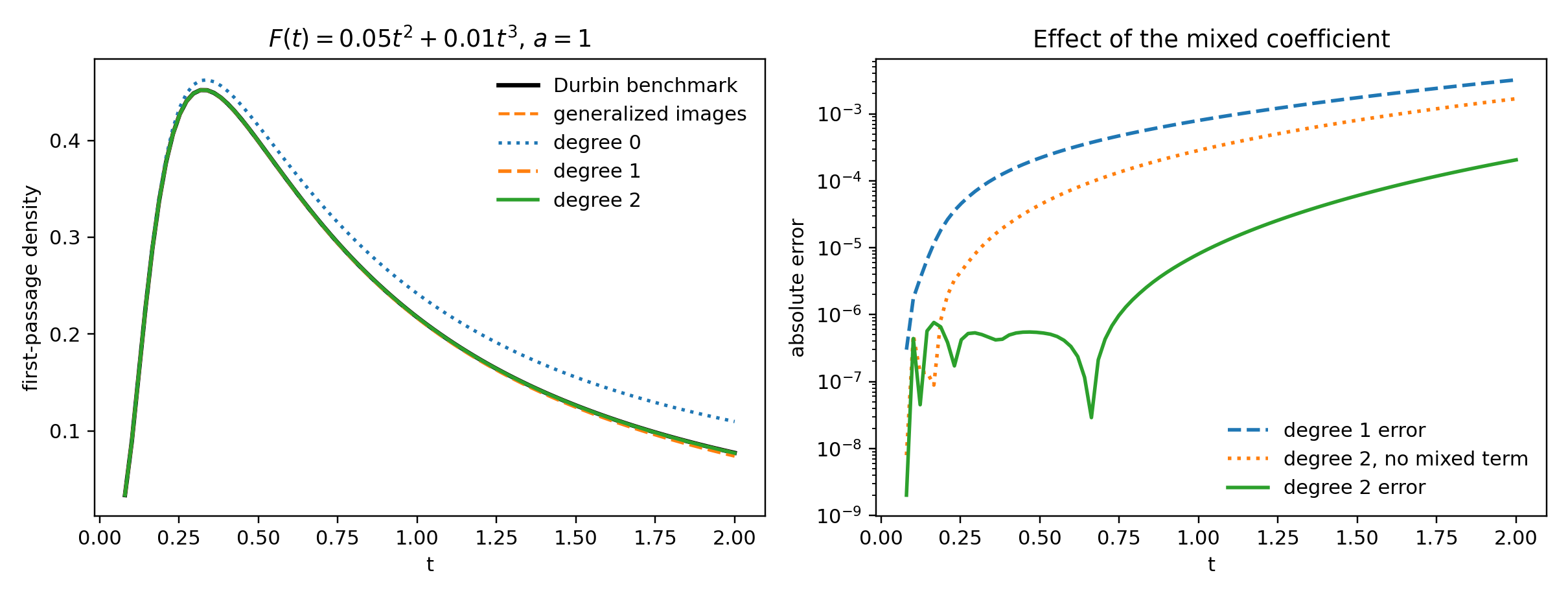}
\caption{Left: the generalized-image flux, the formal transmutation
expansion, and the independent Durbin benchmark.  Right: absolute errors,
showing the contribution of the mixed $ABq_{1,1}$ coefficient.}
\label{fig:numerical}
\end{figure}

As a separate check, truncate \eqref{eq:GammaAB} at $m+n\le6$.  At $t=2$,
direct heat transport gives
\[
 |(P_t\Gamma_{A,B}^{[6]})(F_{A,B}(t))|
 =1.24\times10^{-4},
\]
and
\[
 |\partial_x(P_t\Gamma_{A,B}^{[6]})(F_{A,B}(t))-1|
 =2.92\times10^{-4}.
\]
For $t\le1$, both residuals are below $1.5\times10^{-8}$.  These computations
are numerical validation of the exact formal identities and are not used in
their proofs.

As a test with three active polynomial degrees, we also used
\[
 F(t)=0.03t^2+0.006t^3+0.001t^4,
 \qquad a=1.
\]
With the same image grid, the image flux and an independent 5000-step Durbin
calculation differ by at most $2.63\times10^{-7}$ at
$t=0.2,0.5,1,1.5,2$.  The corresponding boundary residuals are below
$1.5\times10^{-7}$.  This experiment tests the full multi-index kernel rather
than a two-coordinate specialization.

\section{Open problems}
\label{sec:open}

The generalized-image construction supplies the genuine density on a
nonzero interval without requiring convergence of the formal transmutation
series.  Several questions remain.

\begin{enumerate}
\item Is the analytic interval in
Theorem~\ref{thm:polynomial-first-passage} maximal?  In particular, when can
the signed image representation be continued for all times before the
geometric contact time?

\item What is the sharp exponential type of
$n_{\boldsymbol\beta}$?  The proof gives the correct multiscale exponent
$\max_j|\beta_j|^{2/(2j-1)}$, but not its optimal constant.

\item For which $t$, $a$ and $\bc$ does
\[
 \sum_{\bm\ge\mathbf0}\bc^{\bm}q_{\bm}(t;a)
\]
converge to the generalized-image density?  Local analytic dependence of the
Volterra solution suggests a positive answer near $\bc=\mathbf0$, but a
uniform coefficient theorem remains to be written.

\item Can the formal spatial inverse
$\T_{-\bc}^{-1}\delta_a$ be realized as a function, distribution or signed
measure and related directly to the exterior image measure
$\mu_{\bc,a}$?  Proposition~\ref{prop:multi-beta-transform} identifies their
kernel seeds, but does not yet identify the full inverse objects.

\item Outside the domain of ordinary convergence, is the transmutation
expansion asymptotic or Borel summable to the image flux?

\item How does the construction behave for nonmonotone polynomial boundaries,
for which the probability of passage and the terminal term in optional
stopping require additional care?

\item What is the minimal scalar differential order of the mixed
constant-flux function?  Theorem~\ref{thm:mixed-ode} gives an order-five
annihilator, but minimality and its dependence on the set of active boundary
degrees remain open.

\item How should one select an annihilator for a particular application?
The order-five and order-eight equations show that differential order,
parameter valuation, stencil size and transported jet complexity need not be
minimized by the same operator.  A description of the annihilator ideal over
polynomial and formal parameter rings may be more relevant to first-passage
perturbation than a single minimal generator over the rational field.

\item Can the construction be extended to an analytic boundary with
infinitely many nonzero Taylor coefficients?  The finite multi-grading and
finite moment closure used here would then require uniform growth control.
\end{enumerate}

The central distinction is now between existence of the probability density,
which is provided by generalized images, and convergence of its differential
transmutation expansion.  The beta-moment transform shows that these are two
realizations of the same graded structure, but a complete equivalence of the
full inverse objects remains open.

\section{Conclusion}

A finite polynomial moving boundary determines a time-independent signed
exterior image measure.  Its multivariate Wright kernel satisfies an explicit
second-kind Volterra equation, and mixed finite-difference cancellation gives
the exponential control needed to construct the killed Brownian heat kernel
on a nonzero interval.  The boundary flux is therefore an actual
first-passage density, not only a formal expansion.

The same polynomial determines an explicit multi-index spatial
transmutation.  Its spatial grades are $W(\bn)$, its temporal moment orders
are $K(\bn)$, and the balance $W=2K-|\bn|$ makes the differential and temporal
descriptions Laplace-conjugate.  The multivariate beta-moment transform gives
the direct link: it carries the Wright seed of the image kernel to the entire
constant-flux transported-zero function.

Mixed boundary powers therefore produce both mixed image interactions and
mixed differential operators.  The quadratic--cubic example shows that
annihilating ODEs need not be unique: a low-order equation and a low-valuation
equation may both annihilate the distinguished function while having
different solution spaces and behaving differently at the perturbative base
point.  The numerical experiments show
simultaneously that the generalized-image flux agrees with an independent
Durbin calculation and that the first mixed transmutation coefficient is
essential.  The probabilistic and differential theories thus close into a
single finite-polynomial framework.

\appendix
\section{The rational quotient between the two mixed operators}
\label{app:quotient}

For completeness, write
\[
 \mathcal Q_{A,B}=q_3D^3+q_2D^2+q_1D+q_0
\]
and put
\[
 \Delta=2A^5+135A^2B^2x+675B^3.
\]
The exact quotient in \eqref{eq:L8-factor-L5} is given by
\begin{align*}
 q_3&=-\frac5{x^{12}\Delta},\\
 q_2&=\frac{75(10A^5+693A^2B^2x+3375B^3)}{x^{13}\Delta^2},\\
 q_1&=-\frac{6N_1}{x^{14}\Delta^3},\\
 q_0&=\frac{2N_0}{x^{15}\Delta^3},
\end{align*}
where
\begin{align*}
N_1={}&16A^{11}x^3+6500A^{10}+60A^9Bx^2
+2160A^8B^2x^4+900000A^7B^2x\\
&+14850A^6B^3x^3+72900A^5B^4x^5+4387500A^5B^3
+31185000A^4B^4x^2\\
&+729000A^3B^5x^4+303750000A^2B^5x
+1822500AB^6x^3+740390625B^6,
\end{align*}
and
\begin{align*}
N_0={}&1088A^{11}x^3+175500A^{10}+4500A^9Bx^2
+150120A^8B^2x^4+24570000A^7B^2x\\
&+1078650A^6B^3x^3+5175900A^5B^4x^5+118462500A^5B^3
+862650000A^4B^4x^2\\
&+51759000A^3B^5x^4+8292375000A^2B^5x
+129397500AB^6x^3+19990546875B^6.
\end{align*}
Using the Ore composition rule
\[
 (aD^r)\circ(bD^s)
 =\sum_{\ell=0}^r\binom r\ell ab^{(\ell)}D^{r-\ell+s},
\]
direct substitution gives
\[
 \widetilde{\mathcal L}_{A,B}
 -\mathcal Q_{A,B}\circ\mathcal L_{A,B}=0.
\]

\section*{Acknowledgements}
The author thanks colleagues who commented on earlier stages of the
transported-zero and Brownian boundary-crossing program.

\section*{Competing interests}
The author declares no competing interests.

\section*{Data and code availability}
The numerical experiment uses no external data.  The Python code generating
Table~\ref{tab:numerical} and Figure~\ref{fig:numerical}, together with the
exact symbolic scripts verifying the mixed annihilators and Ore
factorization, is available from the author and is supplied with the
manuscript source.

\end{document}